\documentclass[onecolumn,amsmath,amssymb,floatfix]{revtex4}

\usepackage{graphicx}
\usepackage{dcolumn}
\usepackage{bm}
\usepackage{epsfig}

\begin{document}


\title{An iteration based on prime and composite factors}
\author{Kyle Kawagoe}
\email{kkawagoe@uchicago.edu}
\affiliation{Department of Physics, University of California Santa Barbara\\ Santa Barbara, California, USA}
\author{Greg Huber}
\email{huber@kitp.ucsb.edu}
\affiliation{Kavli Institute for Theoretical Physics, and\\Department of Physics, University of California Santa Barbara\\ Santa Barbara, California, USA}

\date{\today}

\begin{abstract}
A discrete map based on the sum of an integer's distinct primes factors and the sum of its other factors is defined and its iteration is studied.  There is a subset of the positive integers that will remain positive under repeated iteration of this map.  Three criteria for these numbers are given.  Numerical evidence is presented which suggests that the growth rate of these surviving numbers takes the form of a non-trivial power-law up to leading order.
\end{abstract}
\maketitle

\section{Introduction}
\par
The prime and composite factors of positive integers allow an interesting iteration to be defined.  By comparing the sum of all the prime factors $(\Pi)$ against the sum of all the other factors $(\mathcal{C})$, the function $X(n)$, given by $$X(n)=\Pi(n)-C(n)+n$$ shifts the input integer $n$ by the difference $\Pi(n)-\mathcal{C}(n).$  This function maps the positive integers to the integers.  Because the number and size of the composite factors can easily dwarf the sum of the prime factors, $X(n)$ is typically negative.  However, in some cases, $\Pi(n)$ and $\mathcal{C}(n)$ can achieve a balance, keeping $X(n)$ greater than zero.  If iterating $X(n)$ results in all iterates being positive, in other words, if \begin{align*}
&X^{(a)}(n)>0\text{, }\forall a>0
\end{align*}
where \[X^{(a)}(n)=
\underbrace{X(X(\dotsb X(n)\dotsb))}_{a\text{ times}}
\]
then we say the input $n$ ``survives."
\section{Survivors}
\par
Here are the first 30 survivors (i.e. surviving inputs) under this iteration of the $X$ function:
$$2,3,5,9,11,21,29,35,43,57,85,123,139,155,161,203,209,221,249,259,265,277,299,323,349,403,411,517,521,553$$
For example, 11 is on this list because,
\begin{align*}
X(11)&=\Pi(11)-\mathcal{C}(11)+11=11-1+11=21\\
X(21)&=\Pi(21)-\mathcal{C}(21)+21=(3+7)-(21+1)+21=9\\
X(9)&=\Pi(9)-\mathcal{C}(9)+9=3-(9+1)+9=2\\
X(2)&=\Pi(2)-\mathcal{C}(2)+2=2-1+2=3\\
X(3)&=\Pi(3)-\mathcal{C}(3)+3=3-1+3=5\\
X(5)&=\Pi(5)-\mathcal{C}(5)+5=5-1+5=9\\
\end{align*}
\par
As shown in Figure 1, the $X$ map sends 2$\rightarrow$3, 3$\rightarrow$5, 5$\rightarrow$9, and 9$\rightarrow$2.  It is interesting that (2 3 5 9) appears to be the only fundamental cycle in the $X$ iteration.  (We take X(1)=0, so 1 is not part of any cycle.)  No other cycle has appeared in our  numerical exploration of the function $X$ up to inputs as large as $n=10^9$.  Additionally, it is conceivable that there is some sequence generated by iterating the $X$ map that goes to infinity, but, again, no such sequence has been found up to $n=10^9$.  We have no proof of these observations, but we can restrict the survivors to a special subset of the positive integers with simple prime factorizations.  In particular, we can prove that one of the following must be true of survivors (a proof is in the Appendix):
\begin{enumerate}
\item[(i)] $n$ is prime.
\item[(ii)] $n=p\cdot q$ where $p$ and $q$ are distinct odd primes.
\item[(iii)] $n=9$
\end{enumerate}

\par However, this is not to say that all numbers satisfying these criteria survive.  For example, not all primes are survivors; take the case n=7:
\begin{align*}
X(7)&=13\\
X(13)&=25\\
X(25)&=4\\
X(4)&=1\\
X(1)&=0\\
\end{align*}

Also, not all biprimes (products of two primes \cite{Conway}) are survivors.  For example, 15 and 33 do not survive.  Therefore, the density of survivors is less than $n/log(n)$ by the prime-number theorem\cite{Ingham}.  In fact, numerical exploration of the survivor sequence indicates their density to be considerably less than $n/log(n)$ (see section IV).

\vspace{.1in}
\par
The above criteria lead to an efficient algorithm and associated tree for generating survivors without resorting to a brute-force search.  Results from this tree-based algorithm agree with results from exhaustive searches up to at least $n=10^6$.  A depiction of this tree is shown in Figure 1.  Note that $1$ is always included as a composite factor $-$ for a prime input $p$, it is the only composite factor.  Therefore, for instance, $ X(p) = p -1 + p = 2p-1$.  Primes directly linked by $X$ iteration (such as 7$\rightarrow$13 above) must necessarlily form finite chains by the following argument.  Let $p_i=X^{(i-1)}(p_1)$.  All such sequences take the form $p_i=2^{i-1}p_1-2^{i-1}+1$.  An application of Fermat's Little Theorem to $p_{p_1}$ gives 

$$\text{ }p_{p_1}=-2^{p_1-1}+1=-1+1=0\text{ }(\text{mod }p),$$

\noindent implying that $p_{p_1}$ is not prime.  This implies that there are no ``prime only" sequences generated by $X$ that go to infinity, and therefore any prime survivors lead to iterations that go through composite numbers.  The finite sequences of primes generated by $X$ have been studied in a different context, and are called ``Cunningham chains of the 2nd kind" \cite{Forbes}.
\par Without going into great detail, we briefly discuss the reasoning behind the three criteria (i)-(iii).  It will be shown in the Appendix that any numbers which have more than three prime factors (not necessarily distinct) must map to a non-positive integer.  This limits our search to primes (i) and biprimes (ii and iii).  We can be even more specific.  All even biprimes will map to a smaller even number.  Therefore, if an even biprime is to survive, it must either map to 2 or a smaller even biprime.  Since the only number that maps to 2 is 9 (an odd number), we can use inductive reasoning to prove that there are no even survivors except for 2.  This further limits the second criterion.  All square biprimes map to an even number.  Therefore, the only biprime with non-distinct prime factors is 9.  Thus, all surviving biprimes must either be odd and have distinct prime factors (ii) or be 9 (iii). 
\begin{figure}
  \centering
  \caption[Two Chain on Larger System]{Flow of the $X$ map represented as a tree.}
  \includegraphics[width=0.4\textwidth]%
    {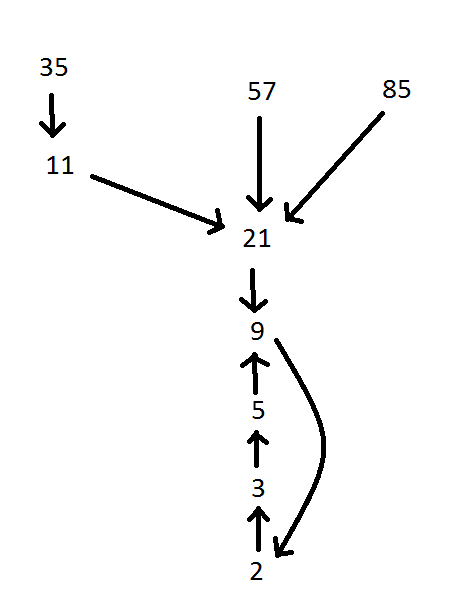}
\end{figure}
\section{Preimage Algorithm}
\par
In order to efficiently and thoroughly search for survivors below some cut off, a preimage algorithm was used.  A quick check shows that the only survivors that map to 2, 3, and 5 are 9, 2, and 3, respectively.  The only survivors that map to 9 are 5 and 21.  All survivors that are not equal to 9 must either be prime or be the product of two distinct odd primes.  Therefore, it is possible to check the preimage of $X$ from any survivor in an informed way.  Let $n$ be a survivor.  If $(n-1)/2$ is prime then $(n-1)/2$ is also a survivor since $X((n-1)/2)=n$.  If $k,j$ are distinct, odd, prime integers such that $n+1=k+j$, then $k\cdot j$ is a survivor since $X(k\cdot j)=n$.  These two cases make up the entirety of $X^{-1}(n)$.

\par
The algorithm to find all survivors between $1$ and $N$ goes as follows.  First, establish a cut off, $N$, that represents the largest possible survivor in the survivor list to be created.  Check all prime numbers $1<p\leq N$ to see if they survive.  If $p$ is found to survive, add it and all of its iterates to the list.  This may include some numbers that are greater than $N$.  Now, for each element in the survivor list that is less than the cut off, namely $n\leq N$, check all primes $2<p<n$ to see if $n-p+1$ is prime and not equal to $p$.  If these conditions are met, add $p\cdot(n-p+1)$ to the list.  Continue this process until all elements of the list less than or equal to $N$ have been checked.  This yields a list of all survivors between 1 and $N$.  This algorithm is illustrated in Figure 2.

\section{Survivor Searching}
\par
If we index the list of survivors $n(k)$ such that 

$$n(1)=2,n(2)=3,n(3)=5,n(4)=9, \ldots $$

\noindent then it is not hard to see that $n(k)>k$. (Since $1\rightarrow0$ under the map $X$, and zero is not a positive integer, 1 is not a survivor.)  Plots of $n(k)/k^{\alpha}$ suggest that the scaling of the survivors is best described by a power-law exponent between $\alpha=1.2$ and $1.5$.  Figure 3 shows a comparison between dividing out different putative power-law behaviors.
\newpage
\section{Conclusion}
\par
The prime and composite factors of the positive integers have been studied for thousands of years.  Iterative algorithms, such as Euclid's Algorithm, operating over the same domain are perhaps just as old.  By combining the integers' multiplicative properties with the operations of addition and subtraction (basically a ``walk" from integer to integer), we have found  a simple map that displays novel scaling behavior.  One can imagine that other nontrivial scalings are lying in wait in this map and in related maps that mix additive and multiplicative properties.

\begin{figure}
  \centering
  \caption[Two Chain on Larger System]{This tree represents the flow of the preimage algorithm.  Notice that 21 branches to 11, 57, and 85.  This is because (21+1)/2=11 is prime, 21$+1=3+19=5+$17, 3$\cdot19=$57 and 5$\cdot 17=$85.}
  \includegraphics[width=0.4\textwidth]%
    {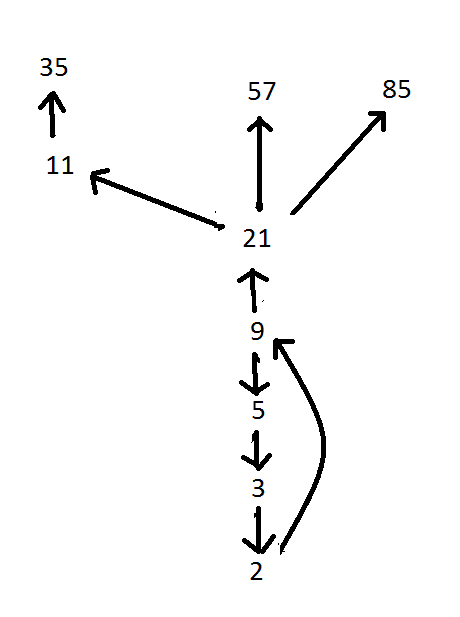}
\end{figure}
\begin{figure}
  \centering
  \caption[Two Chain on Larger System]{This figure shows $n(k)/k^{\alpha}$ against $1/log(k)$ on a log plot for varying values of $\alpha$.   The upper (blue), middle (red), and lower (green) curves correspond to the three exponents  $\alpha =$ 1.5, 1.3, and 1.2, respectively.  This plot suggests that $n(k)=\mathcal{O}(x^{1.5})$ and $x^{1.2}=\mathcal{O}(n(k))$.  We include $\alpha = 1.3$ as a hypothesis for a possible scaling law $n(k)\sim x^{1.3}$.}
  \includegraphics[width=0.4\textwidth]%
    {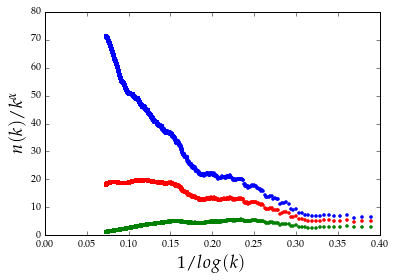}
\end{figure}

\section*{ACKNOWLEDGMENTS}
\noindent This research was supported by the National Science Foundation under grant NSF PHY11-25915.

\clearpage
\appendix

\section{}
\noindent All survivors must meet one of the following criteria.
\begin{enumerate}
\item[(i)] $n$ is prime.
\item[(ii)] $n=p\cdot q$ where $p$ and $q$ are distinct odd primes.
\item[(iii)] $n=9$
\end{enumerate}
\par A series of lemmas will assist in the proof of this statement.  First of all, all positive integers that are the product of $N\geq 3$ distinct prime numbers map to a negative integer (Lemma 1).  If $n=p^3$ where $p$ is prime, then $X(n)=p-(p^3+p^2+1)+p^3=-(p-1)p-1<0$.  Also, if $n=p^2q$ where $p,q$ are distinct primes, then $X(n)<0$ (Lemma 2).  Therefore, if a number is the product of three primes or if it is the product of three or more distinct primes, then it maps to a negative number.  It is also true that if $X(n)\leq 0$ and $p|n$ then $X(pn)<0$ (Lemma 3).  Therefore, by induction, any number that is the product of three or more primes is not a survivor.  This limits the possibilities to primes and biprimes.  By Lemma 4, the only even survivor is 2.  Notice that all square biprimes $p^2$ (except for 4) have the property that $X(p^2)=p-(p^2+1)-p^2=p-1$ is even.  In this exception, $X(X(4))=X(1)=0$ so $4$ does not survive.  Since 2 is the only even survivor, the only square biprime that survives is 9.  This proves that the three criteria are necessary for survivors.\\
Lemma 1: Let $n=\prod\limits_{i=1}^N p_i$ where $N\geq3$ and $p_i$ are distinct primes.  Identify $p_{N+1}=p_1$.
\begin{align*}
\Pi(n)&=\sum\limits_{i=1}^N p_i\\
\mathcal{C}(n)&\geq 1+n+\sum\limits_{i=1}^N p_ip_{i+1}\\
X(n)&=\Pi(n)-\mathcal{C}(n)+n\\
&\leq \sum\limits_{i=1}^Np_i-\left(1+n+\sum\limits_{i=1}^N p_ip_{i+1}\right)+n\\
&\leq 1-\sum\limits_{i=1}^N p_i(p_{i+1}-1)\\
&<0
\end{align*}
Therefore, $n$ is not a survivor.\\\\
Lemma 2:  Let $n=p^aq$ where $p,q$ are distinct primes and $a\geq 2$.
\begin{align*}
X(n)&=p+q-\left(\sum\limits_{i=1}^{a}p^iq+\sum\limits_{i=2}^{a}p^i\right)-p^aq\\
&=-p(p-1)-q(p-1)-\left(\sum\limits_{i=2}^{a-1}p^iq+\sum\limits_{i=3}^{a}p^i\right)\\
&<0
\end{align*}
Lemma 3: Let $n$ be a positive integer such that $X(n)\leq 0$ and $p$ be a prime that divides $n$ exactly $a$ times.  Since $n$ is not a prime, we can rewrite $n=p^a k$ with $a\geq1$.  If $k$ is prime, then by lemma 2, $X(pn)<0$.  Otherwise, $k$ is not prime and $C(k)\geq k$.  In this case, the following calculation can be performed:
\begin{align*}
\mathcal{C}(n)&=\sum\limits_{i=0}^a p^i\mathcal{C}(k)-p\\
\mathcal{C}(pn)&=\sum\limits_{i=0}^{a+1}p^i\mathcal{C}(k)-p=p^{a+1}\mathcal{C}(k)+C(n)\\
\Pi(n)&=\Pi(pn)\\
X(pn)&=\Pi(pn)-\mathcal{C}(pn)+pn\\
&=\Pi(n)-\left(p^{a+1}\mathcal{C}(k)+\mathcal{C}(n)\right)+pn\\
&=X(n)-p^{a+1}(C(k)-k)\\
&<0
\end{align*}
Lemma 4:  From our previous result on the product of three or more primes, if $n$ survives, then $n$ must be prime or biprime. If $n$ is an even biprime, then it takes the form $n=2p$.  Since $n$ is a survivor, $p\neq 2$.  Therefore, $$X(n)=2+p-(2p+1)-2p=p+1$$ which is even and less than $n$, so either $X(n)=2$ or $X(n)=2q$, where $q$ is prime.  The only number that maps to $2$ is $9$, and $9$ is odd.  Therefore, $X(n)=2q$.  Iterating this map, a decreasing sequence containing numbers that are the product of 2 and a prime is created.  However, no such decreasing sequence can have infinite length.  Therefore, $n$ is not a survivor.  It should also be noted that 2 is in the previously mentioned fundamental cycle, so 2 survives.  Thus, 2 is the only even survivor.  This lemma completes the proof of our survivor criteria.

%

\end{document}